\documentclass[11pt,fleqn]{article}

\usepackage[cp1251]{inputenc}

\usepackage{amsmath,amssymb,amsthm,enumerate,epsfig,graphics,cite}
\usepackage[ukrainian,english]{babel}

\newcommand{\todo}[1][\null]{\ensuremath{\clubsuit}}

\newcommand{\noprint}[1]{}

\begin{document}
\begin{center}
\Large\bf
Generalization of Jamet's test for convergence of number series and its new modifications
\end{center}
\begin{center}
\it
Artem M. Ponomarenko
\end{center}

\begin{center}
National Technical University of Ukraine
"Igor Sikorsky Kyiv Polytechnic Institute"
\end{center}

{\it

In this article, we present new generalizations of logarithmic convergence tests for number series, from which we will derive various new generalizations of the Jamet's convergence test. Further, similarly, on the basis of the generalizations of the Schlomilch's test we found, we will obtain modified tests of the convergence of number series.

}
\medskip


\medskip

In \cite {1} the classical Jamet's test is shown, the proof of which can be based on the logarithmic convergence test. Some other tests of convergence of number series, presented in \cite {2}-\cite {5}, are special cases of our generalizations.

In this article we will present generalizations logarithmic test for convergence of number series and the generalized Schlomilch's test.
Next, we show various new generalizations of the Jamet's test and new modified tests for the convergence of number series.

The purpose of our research was to consider and study all types of possible generalizations of criteria for the convergence of number series, which we can consider as a generalization of the classical Jamet's test, as well as all important special cases of these convergence tests.
In view of this, the article will contain many theorems that are similar at first glance, since it is necessary to consider and present each individual case.

\medskip

Throughout the entire article, for any $n \in \mathbb{N}$ we will assume
\begin{gather}\label{eqd1}
\lambda_{0}(n)=1, \quad \lambda_{1}(n)=n, \quad \lambda_{k+1}(n)=\log\left( \lambda_{k}(n) \right), \quad k\in \mathbb{N}.
\end{gather}
Of course, for some first values $n=1$, $2$, ... , $l_{0}$, $l_{0}\in \mathbb{N}$, our sequence $\{\lambda_{k}(n)\mid n \in \mathbb{N}\}$ may not exist for a fixed value of $k$ or may have complex values. However, it is obvious that for any $k$ there exists a number $\gamma\in \mathbb{N}$ such that for any number $n>\gamma$ our sequence $\{\lambda_{k}(n)\mid n \in \mathbb{N}\}$ will exist. Since we are studying the issues of convergence and divergence of number series, it will not affect our research in any way.

Also, let us introduce the notation $\mathbb{N}_{+}:=\mathbb{N}\cup\{ 0 \}$. Here and further, we will assume that $N \in \mathbb{N}$.
Also, here and further, $\delta$, $\delta_{1}$, $\delta_{2}$, ... , $\delta_{7}$ will be some real constant numbers.

Let the functions $f(n)$, $g(n)$, $\varphi(n)$, $\zeta(n)$ be defined for $ \forall n \in \mathbb{N} $ and throughout the article, $n$ will also belong to the set of natural numbers.

We will further assume that some function $\chi(n)$ is positive definite if $\forall n\in \mathbb{N} $
$$
\chi(n)>0;
$$

some function $\chi(n)$ is increasing if $\forall n\in \mathbb{N} $
$$
\chi(n+1)>\chi(n);
$$

some function $\chi(n)$ is unbounded if for $\forall A\in \mathbb{R}$ there exists a number $ M\in \mathbb{N} $ such that
$$
|\chi(M)|\geq A;
$$

some function $\chi(n)$ is non-decreasing if $\forall n\in \mathbb{N} $
$$
\chi(n+1)\geq \chi(n).
$$

\medskip
\section{Generalized tests for the convergence of number series}

First, we need to find generalizations of the test for logarithmic convergence of number series and the Schlomilch's test, which in the future can serve as reference results for proving generalizations of the Jamet's test and modified tests.

{\it Theorem 1 (Generalized Logarithmic test). Let the number sequence} $\{a_{n}\mid n \in \mathbb{N}\}$ {\it is given such that } $a_{n}>0$ {\it for any } $ n \in \mathbb{N}$. {\it Then

1) if there exists a number } $k\in \mathbb{N}_{+}$ {\it such that} $\exists N:$ $\forall n>N$ {\it the following inequality }
\begin{gather}\label{eq1}
L_{k,n}:=\frac{\log \left( \frac{1}{a_{n}\prod\limits_{i=0}^{k}\lambda_{i}(n)}\right)} {\lambda_{k+2}(n)} \geq 1+\delta,  \quad \delta>0
\end{gather}
{\it is satisfied, then the number series } $\sum\limits_{n=1}^{\infty}a_{n}$ {\it converges;

2) if there exists a number } $k\in \mathbb{N}_{+}$ {\it such that} $\exists N:$ $\forall n>N$ {\it the following inequality }
\begin{gather}\label{eq2}
L_{k,n} \leq 1-\delta,  \quad \delta>0
\end{gather}
 {\it is satisfied, then the number series } $\sum\limits_{n=1}^{\infty}a_{n}$ {\it diverges.}

{\it Proof} is based on comparing of the sequence $\{a_{n}\mid n \in \mathbb{N}\}$ with the sequence $\{b_{n,k_{0}}\mid n \in \mathbb{N}\}$, where $b_{n,k_{0}}:=\frac{1}{ \left( \lambda_{k_{0}+1}(n) \right)^{1\pm\delta} \prod\limits_{i=0}^{k_{0}}\lambda_{i}(n)}$, $k_{0}\in \mathbb{N}_{+}$ for $\forall n\in \mathbb{N}$, and further using the logarithm of both sides of the inequality. $\square$

Of course, for some first values $1$, $2$, ... , $l$, $l\in \mathbb{N}$, our sequence $\{b_{n,k_{0}}\mid n \in \mathbb{N}\}$ may not exist due to the definition \eqref{eqd1}, however, obviously, there exists always a sufficiently large number $l_{1}\in \mathbb{N}$ such that for any number $n>l_{1}$, our sequence $\{b_{n,k_{0}}\mid n \in \mathbb{N}\}$ will exist for any $n$ and $k_{0}$. In the future we will not dwell on this.

{\it Theorem 2 (Generalized Schlomilch's test). Let the number sequence} $\{a_{n}\mid n \in \mathbb{N}\}$ {\it is given such that } $a_{n}>0$ {\it for any } $ n \in \mathbb{N}$. {\it Then

1) if there exists a number } $k\in \mathbb{N}_{+}$ {\it such that} $\exists N:$ $\forall n>N$ {\it the following inequality }
\begin{gather}\label{eq3}
S_{k,n}:=\frac{\log \left( \frac{a_{n}}{a_{n+1} } \prod\limits_{i=0}^{k}\frac{\lambda_{i}(n)}{\lambda_{i}(n+1)}\right)} {\lambda_{k+2}(n+1)-\lambda_{k+2}(n)} \geq 1+\delta_{1},  \quad \delta_{1}>0
\end{gather}
{\it is satisfied, then the number series } $\sum\limits_{n=1}^{\infty}a_{n}$ {\it converges;

2) if there exists a number } $k\in \mathbb{N}_{+}$ {\it such that} $\exists N:$ $\forall n>N$ {\it the following inequality }
\begin{gather}\label{eq4}
S_{k,n} \leq 1-\delta_{1},  \quad \delta_{1}>0
\end{gather}
{\it is satisfied, then the number series } $\sum\limits_{n=1}^{\infty}a_{n}$ {\it diverges.}

{\it Proof}  is based on comparing sequence $\{\frac{a_{n+1}}{a_{n}}\mid n \in \mathbb{N}\}$ with sequence $\{\frac{\eta_{n+1,k_{1}}}{\eta_{n,k_{1}}}\mid n \in \mathbb{N}\}$ by variable $n$, where
$$\eta_{n,k_{1}}:=\frac{1}{ \left( \lambda_{k_{1}+1}(n) \right)^{1\pm\delta} \prod\limits_{i=0}^{k_{1}}\lambda_{i}(n)}, \quad k_{1}\in \mathbb{N}_{+}$$
in a similar way. $\square$

{\it Theorem 2$^{\ast}$ (Second version of generalized Schlomilch's test). Let the number sequence} $\quad$ $\{a_{n}\mid n \in \mathbb{N}\}$ {\it is given such that } $a_{n}>0$ {\it for any } $ n \in \mathbb{N}$. {\it Then

1) if there exists a number } $k\in \mathbb{N}_{+}$ {\it such that} $\exists N:$ $\forall n>N$ {\it the following inequality }
\begin{gather}\label{eq3second}
\Omega_{k,n}:=\log \left( \frac{a_{n}}{a_{n+1} } \prod\limits_{i=0}^{k}\frac{\lambda_{i}(n)}{\lambda_{i}(n+1)}\right)\prod\limits_{i=0}^{k+1}\lambda_{i}(n) \geq 1+\delta_{2},  \quad \delta_{2}>0
\end{gather}
{\it is satisfied, then the number series } $\sum\limits_{n=1}^{\infty}a_{n}$ {\it converges;

2) if there exists a number } $k\in \mathbb{N}_{+}$ {\it such that} $\exists N:$ $\forall n>N$ {\it the following inequality }
\begin{gather}\label{eq4second}
\Omega_{k,n} \leq 1-\delta_{2},  \quad \delta_{2}>0
\end{gather}
{\it is satisfied, then the number series } $\sum\limits_{n=1}^{\infty}a_{n}$ {\it diverges.}

The {\it proof} of {\it Theorem 2$^{\ast}$} follows from {\it Theorem 2} and the equivalence of limits $\quad$ $\quad$ $\quad$ $\lambda_{k+1}(n+1)-\lambda_{k+1}(n)$ and $\left(\prod\limits_{i=0}^{k}\lambda_{i}(n)\right)^{-1}$, $k\in \mathbb{N}_{+}$ for  $n\rightarrow +\infty$.   $\square$

For $k=0$, from {\it Theorem 1} we can obtain the classical test for logarithmic convergence, and from {\it Theorem 2} (or {\it Theorem 2$^{\ast}$}) for $k=0$ one can find the classical Schlomilch's test, which is presented in \cite {2}.

Next, we will show generalizations of the Jamet's test and modified tests of convergence for number series.

{\it Theorem 3 (Generalized Jamet's test).} {\it Let the following conditions be satisfied:}

1) {\it the number sequence } $\{a_{n}\mid n \in \mathbb{N}\}$ {\it is given such that } $a_{n}>0$ {\it for any } $ n \in \mathbb{N}$;

2) {\it for any} $n \in \mathbb{N}$ {\it the functions} $f(n)$, $\varphi(n)$ {\it are positive definite, unambiguous, increasing and unbounded;}

3) {\it for any} $n \in \mathbb{N}$ {\it the function} $g(n)$ {\it is positive definite and unambiguous;}

4) {\it for some number} $k\in \mathbb{N}_{+}$ {\it the limit}
\begin{gather}\label{eq5}
\lim\limits_{n\rightarrow \infty} \frac{\varphi(n)g(n)}{f(n)\lambda_{k+2}(n)}=q,   \quad q \in \mathbb{R}\setminus\{0\}
\end{gather}
{\it is satisfied; }
\begin{gather}\label{eq6}
5) \quad \lim\limits_{n\rightarrow \infty} \frac{g(n)}{f(n)}=0.
\end{gather}
 {\it Then

1) if for a given number } $k$ {\it satisfying limit \eqref{eq5},} $\exists N:$ $\forall n>N$ {\it the following inequality }
\begin{gather}\label{eq7}
J_{k,n}:= \left( 1- \left( a_{n}\prod\limits_{i=0}^{k}\lambda_{i}(n) \right)^{\frac{1}{\varphi(n)}} \right)\frac{f(n)}{g(n)} \geq p>\frac{1}{q},   \quad p \in \mathbb{R}
\end{gather}
{\it is satisfied, where $p$ is some constant number, then the number series } $\sum\limits_{n=1}^{\infty}a_{n}$ {\it converges;

2) if for a given number } $k$ {\it satisfying limit \eqref{eq5},} $\exists N:$ $\forall n>N$ {\it the following inequality }
\begin{gather}\label{eq8}
\left( 1- \left( a_{n}\prod\limits_{i=0}^{k}\lambda_{i}(n) \right)^{\frac{1}{\varphi(n)}} \right)\frac{f(n)}{g(n)} \leq p<\frac{1}{q}
\end{gather}
{\it is satisfied, where $p$ is some constant number, then the number series } $\sum\limits_{n=1}^{\infty}a_{n}$ {\it diverges.}

{\it Proof.} First we prove point 1). From the inequality \eqref{eq7} we obtain
$$
a_{n}\prod\limits_{i=0}^{k}\lambda_{i}(n)\leq \left( 1-p\frac{g(n)}{f(n)} \right)^{\varphi(n)},
$$
from where, using elementary actions, we find the inequality
$$
\log \frac{1}{a_{n}\prod\limits_{i=0}^{k}\lambda_{i}(n)}\geq\varphi(n) \log \left( \frac{1}{1-p\frac{g(n)}{f(n)}} \right),
$$
from which further, taking into account the limit \eqref{eq5}, we obtain
$$
\frac{\log \frac{1}{a_{n}\prod\limits_{i=0}^{k}\lambda_{i}(n)}}{\lambda_{k+2}(n)}\geq \frac{\varphi(n)}{\lambda_{k+2}(n)} \log \left( \frac{1}{1-p\frac{g(n)}{f(n)}} \right)=\frac{\varphi(n)}{\lambda_{k+2}(n)}\sum\limits_{j=1}^{\infty}\left( p\frac{g(n)}{f(n)}\right)^{j}=
$$
$$
=\frac{\varphi(n) g(n)}{\lambda_{k+2}(n) f(n)}p+o\left( \frac{\varphi(n) g(n)}{\lambda_{k+2}(n) f(n)}\right) \overset {} {\underset { n \rightarrow \infty }{\longrightarrow} } pq.
$$
In inequality \eqref{eq7} by assumption $p>\frac{1}{q}$, whence it is obvious that there exists a number $\delta_{3}>0$ such that from the last inequality we obtain
$$
\frac{\log \frac{1}{a_{n}\prod\limits_{i=0}^{k}\lambda_{i}(n)}}{\lambda_{k+2}(n)}> 1+\delta_{3}, \quad \delta_{3}>0,
$$
and further from {\it Theorem 1 (Generalized Logarithmic test)}, the number series $\sum\limits_{n=1}^{\infty}a_{n}$  converges.

Next we prove point 2). Similarly we find the inequality
$$
\frac{\log \frac{1}{a_{n}\prod\limits_{i=0}^{k}\lambda_{i}(n)}}{\lambda_{k+2}(n)}\leq \frac{\varphi(n)}{\lambda_{k+2}(n)}\sum\limits_{j=1}^{\infty}\left( p\frac{g(n)}{f(n)}\right)^{j} = \frac{\varphi(n) g(n)}{\lambda_{k+2}(n) f(n)}p +o\left( \frac{\varphi(n) g(n)}{\lambda_{k+2}(n) f(n)}\right) \overset {} {\underset { n \rightarrow \infty }{\longrightarrow} } pq.
$$
In inequality \eqref{eq8} by assumption $p<\frac{1}{q}$, which means that we similarly obtain the inequality
$$
\frac{\log \frac{1}{a_{n}\prod\limits_{i=0}^{k}\lambda_{i}(n)}}{\lambda_{k+2}(n)}< 1-\delta_{3}, \quad \delta_{3}>0,
$$
and further from {\it Theorem 1}, the number series $\sum\limits_{n=1}^{\infty}a_{n}$  diverges. $\square$

Next, obviously, we can consider the case when the number $q$ tends to infinity.

 {\it Theorem 3.1 (Generalized Jamet's test for the case $q \rightarrow +\infty $).} {\it Let the following conditions be satisfied:}

1) {\it the number sequence } $\{a_{n}\mid n \in \mathbb{N}\}$ {\it is given such that } $a_{n}>0$ {\it for any } $ n \in \mathbb{N}$;

2) {\it for any} $n \in \mathbb{N}$ {\it the functions} $f(n)$, $\varphi(n)$ {\it are positive definite, unambiguous, increasing and unbounded;}

3) {\it for any} $n \in \mathbb{N}$ {\it the function} $g(n)$ {\it is positive definite and unambiguous;}

4) {\it for some number} $k\in \mathbb{N}_{+}$ {\it the limit}
\begin{gather}\label{eq5111}
 \lim\limits_{n\rightarrow \infty} \frac{\varphi(n)g(n)}{f(n)\lambda_{k+2}(n)}=+\infty
\end{gather}
{\it is satisfied; }
\begin{gather}\label{eq6111}
5) \quad \lim\limits_{n\rightarrow \infty} \frac{g(n)}{f(n)}=0.
\end{gather}
 {\it Then

1) if for a given number } $k$ {\it satisfying limit \eqref{eq5111},} $\exists N:$ $\forall n>N$ {\it the following inequality }
\begin{gather}\label{eq7111}
J_{k,n}:= \left( 1- \left( a_{n}\prod\limits_{i=0}^{k}\lambda_{i}(n) \right)^{\frac{1}{\varphi(n)}} \right)\frac{f(n)}{g(n)} \geq p> 0
\end{gather}
{\it is satisfied, where $p$ is some constant number, then the number series } $\sum\limits_{n=1}^{\infty}a_{n}$ {\it converges;

2) if for a given number } $k$ {\it satisfying limit \eqref{eq5111},} $\exists N:$ $\forall n>N$ {\it the following inequality }
\begin{gather}\label{eq8111}
\left( 1- \left( a_{n}\prod\limits_{i=0}^{k}\lambda_{i}(n) \right)^{\frac{1}{\varphi(n)}} \right)\frac{f(n)}{g(n)} < 0
\end{gather}
{\it is satisfied, then the number series } $\sum\limits_{n=1}^{\infty}a_{n}$ {\it diverges.}

{\it Proof.} Similar to the previous proof, we consistently find the inequalities
$$
\frac{\log \frac{1}{a_{n}\prod\limits_{i=0}^{k}\lambda_{i}(n)}}{\lambda_{k+2}(n)}> pq  \overset {} {\underset { n \rightarrow \infty }{\longrightarrow} }p\cdot (+\infty)=+\infty>1+\delta_{3}, \quad \delta_{3}>0,
$$
and
$$
\frac{\log \frac{1}{a_{n}\prod\limits_{i=0}^{k}\lambda_{i}(n)}}{\lambda_{k+2}(n)}< -pq  \overset {} {\underset { n \rightarrow \infty }{\longrightarrow} }-p\cdot (+\infty)=-\infty<1-\delta_{3}, \quad \delta_{3}>0.
$$
Then, we also apply the Generalized Logarithmic test. $\square$

Next, we consider the case when $q=0$. In this case, we only have to establish a criterion for the divergence of the number series.

 {\it Theorem 3.2 (Generalized Jamet's test for the case $q=0$).} {\it Let the following conditions be satisfied:}

1) {\it the number sequence } $\{a_{n}\mid n \in \mathbb{N}\}$ {\it is given such that } $a_{n}>0$ {\it for any } $ n \in \mathbb{N}$;

2) {\it for any} $n \in \mathbb{N}$ {\it the functions} $f(n)$, $\varphi(n)$ {\it are positive definite, unambiguous, increasing and unbounded;}

3) {\it for any} $n \in \mathbb{N}$ {\it the function} $g(n)$ {\it is positive definite and unambiguous;}

4) {\it for some number} $k\in \mathbb{N}_{+}$ {\it the limit}
\begin{gather}\label{eq5112}
\lim\limits_{n\rightarrow \infty} \frac{\varphi(n)g(n)}{f(n)\lambda_{k+2}(n)}=0
\end{gather}
{\it is satisfied; }
\begin{gather}\label{eq6112}
5) \quad \lim\limits_{n\rightarrow \infty} \frac{g(n)}{f(n)}=0;
\end{gather}

6) {\it for a given number } $k$ {\it satisfying limit \eqref{eq5112},} $\exists N:$ $\forall n>N$ {\it the following inequality }
\begin{gather}\label{eq8112}
\left( 1- \left( a_{n}\prod\limits_{i=0}^{k}\lambda_{i}(n) \right)^{\frac{1}{\varphi(n)}} \right)\frac{f(n)}{g(n)} \leq p
\end{gather}
{\it is satisfied, where $p$ is some constant number.

Then the number series } $\sum\limits_{n=1}^{\infty}a_{n}$ {\it diverges.

Proof} follows from the inequality
$$
\frac{\log \frac{1}{a_{n}\prod\limits_{i=0}^{k}\lambda_{i}(n)}}{\lambda_{k+2}(n)}< pq  \overset {} {\underset { n \rightarrow \infty }{\longrightarrow} }p\cdot 0=0<1-\delta_{3}, \quad \delta_{3}>0.
$$
$\square$

Now it is obviously easy to see that using similar methods one can obtain a new modified convergence test based on the Schlomilch's test. We will do this also for three cases: $q \in \mathbb{R}\setminus\{0\}$, $q =0$, $q \rightarrow +\infty $.

{\it Theorem 4 (Modification test).} {\it Let the following conditions be satisfied:}

1) {\it the number sequence } $\{a_{n}\mid n \in \mathbb{N}\}$ {\it is given such that } $a_{n}>0$ {\it for any } $ n \in \mathbb{N}$;

2) {\it for any} $n \in \mathbb{N}$ {\it the functions} $f(n)$ {\it is positive definite, unambiguous, increasing and unbounded;}

3) {\it for any} $n \in \mathbb{N}$ {\it the function} $\zeta(n)$ {\it is positive definite, unambiguous and non-decreasing;}

4) {\it for any} $n \in \mathbb{N}$ {\it the function} $g(n)$ {\it is positive definite and unambiguous;}

5) {\it for some number} $k\in \mathbb{N}_{+}$ {\it the limit}
\begin{gather}\label{eq9}
\lim\limits_{n\rightarrow \infty} \frac{\zeta(n)g(n)}{f(n)\left( \lambda_{k+2}(n+1)-\lambda_{k+2}(n)\right) }=s, \quad s \in \mathbb{R}\setminus\{0\}
\end{gather}
{\it is satisfied; }
\begin{gather}\label{eq10}
6) \quad \lim\limits_{n\rightarrow \infty} \frac{g(n)}{f(n)}=0.
\end{gather}
 {\it Then

1) if for a given number } $k$ {\it satisfying limit \eqref{eq9},} $\exists N:$ $\forall n>N$ {\it the following inequality }
\begin{gather}\label{eq11}
P_{k,n}:= \left( 1- \left( \frac{a_{n+1}}{a_{n} } \prod\limits_{i=0}^{k}\frac{\lambda_{i}(n+1)}{\lambda_{i}(n)} \right)^{\frac{1}{\zeta(n)}} \right)\frac{f(n)}{g(n)} \geq r>\frac{1}{s},   \quad r \in \mathbb{R}
\end{gather}
{\it is satisfied, where $r$ is some constant number, then the number series } $\sum\limits_{n=1}^{\infty}a_{n}$ {\it converges;

2) if for a given number } $k$ {\it satisfying limit \eqref{eq9},} $\exists N:$ $\forall n>N$ {\it the following inequality }
\begin{gather}\label{eq12}
\left( 1- \left( \frac{a_{n+1}}{a_{n} } \prod\limits_{i=0}^{k}\frac{\lambda_{i}(n+1)}{\lambda_{i}(n)} \right)^{\frac{1}{\zeta(n)}} \right)\frac{f(n)}{g(n)} \leq r<\frac{1}{s}
\end{gather}
{\it is satisfied, where $r$ is some constant number, then the number series } $\sum\limits_{n=1}^{\infty}a_{n}$ {\it diverges.}

{\it Proof.} First we prove point 1). From the inequality \eqref{eq11} we obtain
$$
\frac{a_{n+1}}{a_{n} } \prod\limits_{i=0}^{k}\frac{\lambda_{i}(n+1)}{\lambda_{i}(n)} \leq \left( 1-r\frac{g(n)}{f(n)} \right)^{\zeta(n)},
$$
from where, using elementary actions, we find the inequality
$$
\log \left( \frac{a_{n}}{a_{n+1}}\prod\limits_{i=0}^{k}\frac{\lambda_{i}(n)}{\lambda_{i}(n+1)} \right) \geq\zeta(n) \log \left( \frac{1}{1-r\frac{g(n)}{f(n)}} \right),
$$
from which further, taking into account the limit \eqref{eq9}, we obtain
$$
\frac{\log \left( \frac{a_{n}}{a_{n+1}}\prod\limits_{i=0}^{k}\frac{\lambda_{i}(n)}{\lambda_{i}(n+1)} \right)}{\lambda_{k+2}(n+1)-\lambda_{k+2}(n)}  \geq \frac{\zeta(n)}{\lambda_{k+2}(n+1)-\lambda_{k+2}(n)} \log \left( \frac{1}{1-r\frac{g(n)}{f(n)}} \right)=
$$
$$
=\frac{\zeta(n)}{ \lambda_{k+2}(n+1)-\lambda_{k+2}(n) }\sum\limits_{j=1}^{\infty}\left( r\frac{g(n)}{f(n)}\right)^{j}=
$$
$$
=\frac{\zeta(n) g(n)}{\left( \lambda_{k+2}(n+1)-\lambda_{k+2}(n) \right) f(n)}r+o\left( \frac{\zeta(n) g(n)}{\left( \lambda_{k+2}(n+1)-\lambda_{k+2}(n)\right) f(n)}\right) \overset {} {\underset { n \rightarrow \infty }{\longrightarrow} } rs.
$$
In inequality \eqref{eq11} by assumption $r>\frac{1}{s}$, whence it is obvious that there exists a number $\delta_{4}>0$ such that from the last inequality we obtain
$$
\frac{\log \left( \frac{a_{n}}{a_{n+1} } \prod\limits_{i=0}^{k}\frac{\lambda_{i}(n)}{\lambda_{i}(n+1)}\right)} {\lambda_{k+2}(n+1)-\lambda_{k+2}(n)}> 1+\delta_{4}, \quad \delta_{4}>0,
$$
and further from {\it Theorem 2 (Generalized Schlomilch's test)}, the number series $\sum\limits_{n=1}^{\infty}a_{n}$  converges.

Next we prove point 2). Similarly we find the inequality
$$
\frac{\log \left( \frac{a_{n}}{a_{n+1} } \prod\limits_{i=0}^{k}\frac{\lambda_{i}(n)}{\lambda_{i}(n+1)}\right)} {\lambda_{k+2}(n+1)-\lambda_{k+2}(n)} \leq \frac{\zeta(n)}{\lambda_{k+2}(n)}\sum\limits_{j=1}^{\infty}\left( r\frac{g(n)}{f(n)}\right)^{j}=
$$
$$
\frac{\zeta(n) g(n)}{\left( \lambda_{k+2}(n+1)-\lambda_{k+2}(n)\right) f(n)}r +o\left( \frac{\zeta(n) g(n)}{\left( \lambda_{k+2}(n+1)-\lambda_{k+2}(n)\right) f(n)}\right) \overset {} {\underset { n \rightarrow \infty }{\longrightarrow} } rs.
$$
In inequality \eqref{eq12} by assumption $r<\frac{1}{s}$, which means that we similarly obtain the inequality
$$
\frac{\log \left( \frac{a_{n}}{a_{n+1} } \prod\limits_{i=0}^{k}\frac{\lambda_{i}(n)}{\lambda_{i}(n+1)}\right)} {\lambda_{k+2}(n+1)-\lambda_{k+2}(n)}< 1-\delta_{4}, \quad \delta_{4}>0,
$$
and from {\it Theorem 2}, the number series $\sum\limits_{n=1}^{\infty}a_{n}$  diverges. $\square$

{\it Theorem 4.1 (Modification test for the case} $ s\rightarrow +\infty $ {\it ).  Let the following conditions be satisfied:}

1) {\it the number sequence } $\{a_{n}\mid n \in \mathbb{N}\}$ {\it is given such that } $a_{n}>0$ {\it for any } $ n \in \mathbb{N}$;

2) {\it for any} $n \in \mathbb{N}$ {\it the functions} $f(n)$ {\it is positive definite, unambiguous, increasing and unbounded;}

3) {\it for any} $n \in \mathbb{N}$ {\it the function} $\zeta(n)$ {\it is positive definite, unambiguous and non-decreasing;}

4) {\it for any} $n \in \mathbb{N}$ {\it the function} $g(n)$ {\it is positive definite and unambiguous;}

5) {\it for some number} $k\in \mathbb{N}_{+}$ {\it the limit}
\begin{gather}\label{eq9111}
\lim\limits_{n\rightarrow \infty} \frac{\zeta(n)g(n)}{f(n)\left( \lambda_{k+2}(n+1)-\lambda_{k+2}(n)\right) }=+\infty
\end{gather}
{\it is satisfied; }
\begin{gather}\label{eq10111}
6) \quad \lim\limits_{n\rightarrow \infty} \frac{g(n)}{f(n)}=0.
\end{gather}
 {\it Then

1) if for a given number } $k$ {\it satisfying limit \eqref{eq9111},} $\exists N:$ $\forall n>N$ {\it the following inequality }
\begin{gather}\label{eq11111}
P_{k,n}:= \left( 1- \left( \frac{a_{n+1}}{a_{n} } \prod\limits_{i=0}^{k}\frac{\lambda_{i}(n+1)}{\lambda_{i}(n)} \right)^{\frac{1}{\zeta(n)}} \right)\frac{f(n)}{g(n)} \geq r>0
\end{gather}
{\it is satisfied, where $r$ is some constant number, then the number series } $\sum\limits_{n=1}^{\infty}a_{n}$ {\it converges;

2) if for a given number } $k$ {\it satisfying limit \eqref{eq9111},} $\exists N:$ $\forall n>N$ {\it the following inequality }
\begin{gather}\label{eq12111}
\left( 1- \left( \frac{a_{n+1}}{a_{n} } \prod\limits_{i=0}^{k}\frac{\lambda_{i}(n+1)}{\lambda_{i}(n)} \right)^{\frac{1}{\zeta(n)}} \right)\frac{f(n)}{g(n)} <0
\end{gather}
{\it is satisfied, then the number series } $\sum\limits_{n=1}^{\infty}a_{n}$ {\it diverges.}

{\it Proof.} Similar to the previous proof, we consistently find the inequalities
$$
\frac{\log \left( \frac{a_{n}}{a_{n+1} } \prod\limits_{i=0}^{k}\frac{\lambda_{i}(n)}{\lambda_{i}(n+1)}\right)} {\lambda_{k+2}(n+1)-\lambda_{k+2}(n)}> rs  \overset {} {\underset { n \rightarrow \infty }{\longrightarrow} }r\cdot (+\infty)=+\infty>1+\delta_{4}, \quad \delta_{4}>0,
$$
and
$$
\frac{\log \left( \frac{a_{n}}{a_{n+1} } \prod\limits_{i=0}^{k}\frac{\lambda_{i}(n)}{\lambda_{i}(n+1)}\right)} {\lambda_{k+2}(n+1)-\lambda_{k+2}(n)}< -rs  \overset {} {\underset { n \rightarrow \infty }{\longrightarrow} }-r\cdot (+\infty)=-\infty<1-\delta_{4}, \quad \delta_{4}>0.
$$
Then, we also apply the Generalized Schlomilch's test. $\square$

{\it Theorem 4.2 (Modification test for the case} $ s\rightarrow 0 $ {\it ).  Let the following conditions be satisfied:}

1) {\it the number sequence } $\{a_{n}\mid n \in \mathbb{N}\}$ {\it is given such that } $a_{n}>0$ {\it for any } $ n \in \mathbb{N}$;

2) {\it for any} $n \in \mathbb{N}$ {\it the functions} $f(n)$ {\it is positive definite, unambiguous, increasing and unbounded;}

3) {\it for any} $n \in \mathbb{N}$ {\it the function} $\zeta(n)$ {\it is positive definite, unambiguous and non-decreasing;}

4) {\it for any} $n \in \mathbb{N}$ {\it the function} $g(n)$ {\it is positive definite and unambiguous;}

5) {\it for some number} $k\in \mathbb{N}_{+}$ {\it the limit}
\begin{gather}\label{eq9112}
\lim\limits_{n\rightarrow \infty} \frac{\zeta(n)g(n)}{f(n)\left( \lambda_{k+2}(n+1)-\lambda_{k+2}(n)\right) }=0
\end{gather}
{\it is satisfied; }
\begin{gather}\label{eq10112}
6) \quad \lim\limits_{n\rightarrow \infty} \frac{g(n)}{f(n)}=0;
\end{gather}
7) {\it for a given number } $k$ {\it satisfying limit \eqref{eq9112},} $\exists N:$ $\forall n>N$ {\it the following inequality }
\begin{gather}\label{eq12112}
\left( 1- \left( \frac{a_{n+1}}{a_{n} } \prod\limits_{i=0}^{k}\frac{\lambda_{i}(n+1)}{\lambda_{i}(n)} \right)^{\frac{1}{\zeta(n)}} \right)\frac{f(n)}{g(n)} <r
\end{gather}
{\it is satisfied, where $r$ is some constant number.

Then the number series } $\sum\limits_{n=1}^{\infty}a_{n}$ {\it diverges.

Proof} follows from the inequality
$$
\frac{\log \left( \frac{a_{n}}{a_{n+1} } \prod\limits_{i=0}^{k}\frac{\lambda_{i}(n)}{\lambda_{i}(n+1)}\right)} {\lambda_{k+2}(n+1)-\lambda_{k+2}(n)}< rs  \overset {} {\underset { n \rightarrow \infty }{\longrightarrow} }r\cdot 0=0<1-\delta_{4}, \quad \delta_{4}>0.
$$
$\square$

In the practical study of series, the special case of the modified convergence test will be most widespread, when $\zeta(n)\equiv 1$.
Let us present all three similar theorems for $\zeta(n)\equiv 1$, which, as we see, do not need proof.

{\it Theorem 5 (Modification test).} {\it Let the following conditions be satisfied:}

1) {\it the number sequence } $\{a_{n}\mid n \in \mathbb{N}\}$ {\it is given such that } $a_{n}>0$ {\it for any } $ n \in \mathbb{N}$;

2) {\it for any} $n \in \mathbb{N}$ {\it the functions} $f(n)$ {\it is positive definite, unambiguous, increasing and unbounded;}

3) {\it for any} $n \in \mathbb{N}$ {\it the function} $g(n)$ {\it is positive definite and unambiguous;}

4) {\it for some number} $k\in \mathbb{N}_{+}$ {\it the limit}
\begin{gather}\label{er9}
\lim\limits_{n\rightarrow \infty} \frac{g(n)}{f(n)\left( \lambda_{k+2}(n+1)-\lambda_{k+2}(n) \right) }=\beta, \quad \beta \in \mathbb{R}\setminus\{0\}
\end{gather}
{\it is satisfied; }
\begin{gather}\label{er10}
5) \quad \lim\limits_{n\rightarrow \infty} \frac{g(n)}{f(n)}=0.
\end{gather}
 {\it Then

1) if for a given number } $k$ {\it satisfying limit \eqref{er9},} $\exists N:$ $\forall n>N$ {\it the following inequality }
\begin{gather}\label{er11}
\tilde{P}_{k,n}:= \left( 1- \frac{a_{n+1}}{a_{n} } \prod\limits_{i=0}^{k}\frac{\lambda_{i}(n+1)}{\lambda_{i}(n)}  \right)\frac{f(n)}{g(n)} \geq \alpha>\frac{1}{\beta},    \quad \alpha \in \mathbb{R}
\end{gather}
{\it is satisfied, where $\alpha$ is some constant number, then the number series } $\sum\limits_{n=1}^{\infty}a_{n}$ {\it converges;

2) if for a given number } $k$ {\it satisfying limit \eqref{er9},} $\exists N:$ $\forall n>N$ {\it the following inequality }
\begin{gather}\label{er12}
\left( 1-  \frac{a_{n+1}}{a_{n} } \prod\limits_{i=0}^{k}\frac{\lambda_{i}(n+1)}{\lambda_{i}(n)}  \right)\frac{f(n)}{g(n)} \leq \alpha<\frac{1}{\beta}
\end{gather}
{\it is satisfied, where $\alpha$ is some constant number, then the number series } $\sum\limits_{n=1}^{\infty}a_{n}$ {\it diverges.}

 {\it Theorem 5.1 (Modification test for the case} $ \beta\rightarrow +\infty ${\it ).  Let the following conditions be satisfied:}

1) {\it the number sequence } $\{a_{n}\mid n \in \mathbb{N}\}$ {\it is given such that } $a_{n}>0$ {\it for any } $ n \in \mathbb{N}$;

2) {\it for any} $n \in \mathbb{N}$ {\it the functions} $f(n)$ {\it is positive definite, unambiguous, increasing and unbounded;}

3) {\it for any} $n \in \mathbb{N}$ {\it the function} $g(n)$ {\it is positive definite and unambiguous;}

4) {\it for some number} $k\in \mathbb{N}_{+}$ {\it the limit}
\begin{gather}\label{er9111}
\lim\limits_{n\rightarrow \infty} \frac{g(n)}{f(n)\left( \lambda_{k+2}(n+1)-\lambda_{k+2}(n)\right) }=+\infty
\end{gather}
{\it is satisfied; }
\begin{gather}\label{er10111}
5) \quad \lim\limits_{n\rightarrow \infty} \frac{g(n)}{f(n)}=0.
\end{gather}
 {\it Then

1) if for a given number } $k$ {\it satisfying limit \eqref{er9111},} $\exists N:$ $\forall n>N$ {\it the following inequality }
\begin{gather}\label{er11111}
\tilde{P}_{k,n}:= \left( 1-  \frac{a_{n+1}}{a_{n} } \prod\limits_{i=0}^{k}\frac{\lambda_{i}(n+1)}{\lambda_{i}(n)} \right)\frac{f(n)}{g(n)} \geq \alpha>0
\end{gather}
{\it is satisfied, where $\alpha$ is some constant number, then the number series } $\sum\limits_{n=1}^{\infty}a_{n}$ {\it converges;

2) if for a given number } $k$ {\it satisfying limit \eqref{er9111},} $\exists N:$ $\forall n>N$ {\it the following inequality }
\begin{gather}\label{er12111}
\left( 1- \frac{a_{n+1}}{a_{n} } \prod\limits_{i=0}^{k}\frac{\lambda_{i}(n+1)}{\lambda_{i}(n)}  \right)\frac{f(n)}{g(n)} <0
\end{gather}
{\it is satisfied, then the number series } $\sum\limits_{n=1}^{\infty}a_{n}$ {\it diverges.}

{\it Theorem 5.2 (Modification test for the case} $ \beta\rightarrow 0 ${\it ).  Let the following conditions be satisfied:}

1) {\it the number sequence } $\{a_{n}\mid n \in \mathbb{N}\}$ {\it is given such that } $a_{n}>0$ {\it for any } $ n \in \mathbb{N}$;

2) {\it for any} $n \in \mathbb{N}$ {\it the functions} $f(n)$ {\it is positive definite, unambiguous, increasing and unbounded;}

3) {\it for any} $n \in \mathbb{N}$ {\it the function} $g(n)$ {\it is positive definite and unambiguous;}

4) {\it for some number} $k\in \mathbb{N}_{+}$ {\it the limit}
\begin{gather}\label{er9112}
\lim\limits_{n\rightarrow \infty} \frac{g(n)}{f(n)\left( \lambda_{k+2}(n+1)-\lambda_{k+2}(n)\right) }=0
\end{gather}
{\it is satisfied; }
\begin{gather}\label{er10112}
5) \quad \lim\limits_{n\rightarrow \infty} \frac{g(n)}{f(n)}=0,
\end{gather}
6) {\it for a given number } $k$ {\it satisfying limit \eqref{er9112},} $\exists N:$ $\forall n>N$ {\it the following inequality }
\begin{gather}\label{er12112}
\left( 1-  \frac{a_{n+1}}{a_{n} } \prod\limits_{i=0}^{k}\frac{\lambda_{i}(n+1)}{\lambda_{i}(n)}  \right)\frac{f(n)}{g(n)} <\alpha
\end{gather}
{\it is satisfied, where $\alpha$ is some constant number.

Then the number series } $\sum\limits_{n=1}^{\infty}a_{n}$ {\it diverges.}

From Theorem 5 for $f(n)=n$, $g(n)=1$, $k=0$ we obtain the classical version of Raabe's test.
From Theorem 3 for $f(n)=n$, $g(n)=1$, $\varphi(n)=n$, $k=0$ we obtain the classical Jamet's test.

{\it Remark 1.} 1) Condition \eqref{eq9} in the {\it Theorem 4} can be replaced by condition
$$
\lim\limits_{n\rightarrow \infty} \frac{g(n)\zeta(n)}{f(n) }\prod\limits_{i=0}^{k+1}\lambda_{i}(n)=s, \quad s \in \mathbb{R}\setminus\{0\};
$$
2) Condition \eqref{eq9111} in the {\it Theorem 4.1} can be replaced by condition $\lim\limits_{n\rightarrow \infty} \frac{g(n)\zeta(n)}{f(n) }\prod\limits_{i=0}^{k+1}\lambda_{i}(n)=+\infty$; $\quad$
3) Condition \eqref{eq9112} in the {\it Theorem 4.2} can be replaced by condition $\lim\limits_{n\rightarrow \infty} \frac{g(n)\zeta(n)}{f(n) }\prod\limits_{i=0}^{k+1}\lambda_{i}(n)=0$; $\quad$
4) Condition \eqref{er9} in {\it Theorem 5} can be replaced by condition
$$
\lim\limits_{n\rightarrow \infty} \frac{g(n)}{f(n) }\prod\limits_{i=0}^{k+1}\lambda_{i}(n)=\beta, \quad \beta \in \mathbb{R}\setminus\{0\};
$$
5) Condition \eqref{er9111} in the {\it Theorem 5.1} can be replaced by condition $\lim\limits_{n\rightarrow \infty} \frac{g(n)}{f(n) }\prod\limits_{i=0}^{k+1}\lambda_{i}(n)=+\infty$; $\quad$ $\quad$
6) Condition \eqref{er9112} in the {\it Theorem 5.2} can be replaced by condition $\lim\limits_{n\rightarrow \infty} \frac{g(n)}{f(n) }\prod\limits_{i=0}^{k+1}\lambda_{i}(n)=0$.

\medskip

\medskip
\section{More general convergence tests for number series}

This section of the article will present further generalizations of all the theorems of the first section. Their proofs are completely similar to those of the previous tests for the convergence of number series shown in Section 1.

 {\it Definition 1.} The function $\psi(x)$ belongs to class $\Phi$ if there exists a number $ \sigma >0$ such that for any $x \in [\sigma,+\infty)$ the following conditions are satisfied:

  1) the function $\psi(x)$ is positive definite and unambiguous;

  2) the function $\psi'(x)$ is continuous;

  3) $\psi'(x)>0$;

  4) $\lim\limits_{x\rightarrow +\infty}  \psi(x) = +\infty $.

  We will denote this as: $\psi \in \Phi$.

Throughout this article, for any function $\psi \in \Phi$ we will assume
\begin{gather}\label{eqd1de}
\lambda_{0}(\psi(n))=1, \quad \lambda_{1}\left(\psi(n)\right)=\psi(n), \quad \lambda_{k+1}\left(\psi(n)\right)=\log\left( \lambda_{k}\left(\psi(n)\right) \right), \quad k\in \mathbb{N}.
\end{gather}

{\it Theorem 6 (Generalized Logarithmic test). Let the number sequence} $\{a_{n}\mid n \in \mathbb{N}\}$ {\it is given such that } $a_{n}>0$ {\it for any } $ n \in \mathbb{N}$ {\it and } $\psi \in \Phi$.  {\it Then

1) if there exists a number } $k\in \mathbb{N}_{+}$ {\it such that} $\exists N:$ $\forall n>N$ {\it the following inequality }
\begin{gather}\label{eq1d}
L_{k,n;\psi}:=\frac{\log \left( \frac{\psi'(n)}{a_{n}\prod\limits_{i=0}^{k}\lambda_{i}\left(\psi(n)\right)}\right)} {\lambda_{k+2}\left(\psi(n)\right)} \geq 1+\delta_{5},  \quad \delta_{5}>0
\end{gather}
{\it is satisfied, then the number series } $\sum\limits_{n=1}^{\infty}a_{n}$ {\it converges;

2) if there exists a number } $k\in \mathbb{N}_{+}$ {\it such that} $\exists N:$ $\forall n>N$ {\it the following inequality }
\begin{gather}\label{eq2d}
L_{k,n;\psi} \leq 1-\delta_{5},  \quad \delta_{5}>0
\end{gather}
 {\it is satisfied, then the number series } $\sum\limits_{n=1}^{\infty}a_{n}$ {\it diverges.}

{\it Theorem 7 (Generalized Schlemilch's test). Let the number sequence} $\{a_{n}\mid n \in \mathbb{N}\}$ {\it is given such that } $a_{n}>0$ {\it for any } $ n \in \mathbb{N}$ {\it and } $\psi \in \Phi$.  {\it Then

1) if there exists a number } $k\in \mathbb{N}_{+}$ {\it such that} $\exists N:$ $\forall n>N$ {\it the following inequality }
\begin{gather}\label{eq3d}
S_{k,n;\psi}:=\frac{\log \left( \frac{a_{n}}{a_{n+1}} \frac{\psi'(n+1)}{\psi'(n) } \prod\limits_{i=0}^{k}\frac{\lambda_{i}\left(\psi(n)\right)}{\lambda_{i}\left(\psi(n+1)\right)}\right)} {\lambda_{k+2}\left(\psi(n+1)\right)-\lambda_{k+2}\left(\psi(n)\right)} \geq 1+\delta_{6},  \quad \delta_{6}>0
\end{gather}
{\it is satisfied, then the number series } $\sum\limits_{n=1}^{\infty}a_{n}$ {\it converges;

2) if there exists a number } $k\in \mathbb{N}_{+}$ {\it such that} $\exists N:$ $\forall n>N$ {\it the following inequality }
\begin{gather}\label{eq4d}
S_{k,n;\psi} \leq 1-\delta_{6},  \quad \delta_{6}>0
\end{gather}
{\it is satisfied, then the number series } $\sum\limits_{n=1}^{\infty}a_{n}$ {\it diverges.}

{\it Theorem 7$^{\ast}$ (Second version of generalized Schlomilch's test). Let the number sequence} $\quad$ $\{a_{n}\mid n \in \mathbb{N}\}$ {\it is given such that } $a_{n}>0$ {\it for any } $ n \in \mathbb{N}$ {\it and}

$\quad $$ 1) \quad \psi \in \Phi;$
\begin{gather}\label{eq3nnfgnn01}
2) \quad \lim\limits_{n\rightarrow \infty} \frac{\psi(n)}{\psi(n+1)-\psi(n)}=  +\infty.
\end{gather}
{\it Then

1) if there exists a number } $k\in \mathbb{N}_{+}$ {\it such that} $\exists N:$ $\forall n>N$ {\it the following inequality }
\begin{gather}\label{eq3second11}
\Omega_{k,n;\psi}:=\log \left( \frac{a_{n}}{a_{n+1} }\frac{\psi'(n+1)}{\psi'(n)} \prod\limits_{i=0}^{k}\frac{\lambda_{i}(\psi(n))}{\lambda_{i}(\psi(n+1))}\right)\frac{\prod\limits_{i=0}^{k+1}\lambda_{i}(\psi(n))}{\psi(n+1)-\psi(n)} \geq 1+\delta_{7},  \quad \delta_{7}>0
\end{gather}
{\it is satisfied, then the number series } $\sum\limits_{n=1}^{\infty}a_{n}$ {\it converges;

2) if there exists a number } $k\in \mathbb{N}_{+}$ {\it such that} $\exists N:$ $\forall n>N$ {\it the following inequality }
\begin{gather}\label{eq4second11}
\Omega_{k,n;\psi} \leq 1-\delta_{7},  \quad \delta_{7}>0
\end{gather}
{\it is satisfied, then the number series } $\sum\limits_{n=1}^{\infty}a_{n}$ {\it diverges.}

{\it Theorem 8 (Generalized Jamet's test).} {\it Let the following conditions be satisfied:}

1) {\it the number sequence } $\{a_{n}\mid n \in \mathbb{N}\}$ {\it is given such that } $a_{n}>0$ {\it for any } $ n \in \mathbb{N}$;

2) {\it for any} $n \in \mathbb{N}$ {\it the functions} $f(n)$, $\varphi(n)$ {\it are positive definite, unambiguous, increasing and unbounded;}

3) {\it for any} $n \in \mathbb{N}$ {\it the function} $g(n)$ {\it is positive definite and unambiguous;}

4) $\psi \in \Phi$;

5) {\it for some number} $k\in \mathbb{N}_{+}$ {\it the limit}
\begin{gather}\label{eq21}
\lim\limits_{n\rightarrow \infty} \frac{\varphi(n)g(n)}{f(n)\lambda_{k+2}(\psi(n))}=Q, \quad Q \in \mathbb{R}\setminus\{0\}
\end{gather}
{\it is satisfied; }
\begin{gather}\label{eq22}
6) \quad \lim\limits_{n\rightarrow \infty} \frac{g(n)}{f(n)}=0.
\end{gather}
 {\it Then

1) if for a given number } $k$ {\it satisfying limit \eqref{eq21},} $\exists N:$ $\forall n>N$ {\it the following inequality }
\begin{gather}\label{eq23}
J_{k,n;\psi}:= \left( 1- \left( \frac{a_{n}}{\psi'(n)}\prod\limits_{i=0}^{k}\lambda_{i}(\psi(n)) \right)^{\frac{1}{\varphi(n)}} \right)\frac{f(n)}{g(n)} \geq T>\frac{1}{Q},    \quad T \in \mathbb{R}
\end{gather}
{\it is satisfied, where $T$ is some constant number, then the number series } $\sum\limits_{n=1}^{\infty}a_{n}$ {\it converges;

2) if for a given number } $k$ {\it satisfying limit \eqref{eq21},} $\exists N:$ $\forall n>N$ {\it the following inequality }
\begin{gather}\label{eq24}
\left( 1- \left( \frac{a_{n}}{\psi'(n)}\prod\limits_{i=0}^{k}\lambda_{i}(\psi(n)) \right)^{\frac{1}{\varphi (n)}} \right)\frac{f(n)}{g(n)} \leq T<\frac{1}{Q}
\end{gather}
{\it is satisfied, where $T$ is some constant number, then the number series } $\sum\limits_{n=1}^{\infty}a_{n}$ {\it diverges.}

{\it Theorem 8.1 (Generalized Jamet's test for the case} $ Q\rightarrow +\infty ${\it ).  Let the following conditions be satisfied:}

1) {\it the number sequence } $\{a_{n}\mid n \in \mathbb{N}\}$ {\it is given such that } $a_{n}>0$ {\it for any } $ n \in \mathbb{N}$;

2) {\it for any} $n \in \mathbb{N}$ {\it the functions} $f(n)$, $\varphi(n)$ {\it are positive definite, unambiguous, increasing and unbounded;}

3) {\it for any} $n \in \mathbb{N}$ {\it the function} $g(n)$ {\it is positive definite and unambiguous;}

4) $\psi \in \Phi$;

5) {\it for some number} $k\in \mathbb{N}_{+}$ {\it the limit}
\begin{gather}\label{eq21111}
\lim\limits_{n\rightarrow \infty} \frac{\varphi(n)g(n)}{f(n)\lambda_{k+2}(\psi(n))}=+\infty
\end{gather}
{\it is satisfied; }
\begin{gather}\label{eq22111}
6) \quad \lim\limits_{n\rightarrow \infty} \frac{g(n)}{f(n)}=0.
\end{gather}
 {\it Then

1) if for a given number } $k$ {\it satisfying limit \eqref{eq21111},} $\exists N:$ $\forall n>N$ {\it the following inequality }
\begin{gather}\label{eq23111}
J_{k,n;\psi}:= \left( 1- \left( \frac{a_{n}}{\psi'(n)}\prod\limits_{i=0}^{k}\lambda_{i}(\psi(n)) \right)^{\frac{1}{\varphi(n)}} \right)\frac{f(n)}{g(n)} \geq T>0
\end{gather}
{\it is satisfied, where $T$ is some constant number, then the number series } $\sum\limits_{n=1}^{\infty}a_{n}$ {\it converges;

2) if for a given number } $k$ {\it satisfying limit \eqref{eq21111},} $\exists N:$ $\forall n>N$ {\it the following inequality }
\begin{gather}\label{eq24111}
\left( 1- \left( \frac{a_{n}}{\psi'(n)}\prod\limits_{i=0}^{k}\lambda_{i}(\psi(n)) \right)^{\frac{1}{\varphi (n)}} \right)\frac{f(n)}{g(n)} <0
\end{gather}
{\it is satisfied, then the number series } $\sum\limits_{n=1}^{\infty}a_{n}$ {\it diverges.}

{\it Theorem 8.2 (Generalized Jamet's test for the case} $ Q\rightarrow 0 ${\it ).  Let the following conditions be satisfied:}

1) {\it the number sequence } $\{a_{n}\mid n \in \mathbb{N}\}$ {\it is given such that } $a_{n}>0$ {\it for any } $ n \in \mathbb{N}$;

2) {\it for any} $n \in \mathbb{N}$ {\it the functions} $f(n)$, $\varphi(n)$ {\it are positive definite, unambiguous, increasing and unbounded;}

3) {\it for any} $n \in \mathbb{N}$ {\it the function} $g(n)$ {\it is positive definite and unambiguous;}

4) $\psi \in \Phi$;

5) {\it for some number} $k\in \mathbb{N}_{+}$ {\it the limit}
\begin{gather}\label{eq21112}
\lim\limits_{n\rightarrow \infty} \frac{\varphi(n)g(n)}{f(n)\lambda_{k+2}(\psi(n))}=0
\end{gather}
{\it is satisfied; }
\begin{gather}\label{eq22112}
6) \quad \lim\limits_{n\rightarrow \infty} \frac{g(n)}{f(n)}=0;
\end{gather}
7) {\it for a given number } $k$ {\it satisfying limit \eqref{eq21112},} $\exists N:$ $\forall n>N$ {\it the following inequality }
\begin{gather}\label{eq24112}
\left( 1- \left( \frac{a_{n}}{\psi'(n)}\prod\limits_{i=0}^{k}\lambda_{i}(\psi(n)) \right)^{\frac{1}{\varphi (n)}} \right)\frac{f(n)}{g(n)} \leq T
\end{gather}
{\it is satisfied, where $T$ is some constant number.

Then the number series } $\sum\limits_{n=1}^{\infty}a_{n}$ {\it diverges.}

{\it Theorem 9 (Modification test).} {\it Let the following conditions be satisfied:}

1) {\it the number sequence } $\{a_{n}\mid n \in \mathbb{N}\}$ {\it is given such that } $a_{n}>0$ {\it for any } $ n \in \mathbb{N}$;

2) {\it for any} $n \in \mathbb{N}$ {\it the function} $f(n)$ {\it is positive definite, unambiguous, increasing and unbounded;}

3) {\it for any} $n \in \mathbb{N}$ {\it the function} $\zeta(n)$ {\it is positive definite, unambiguous and non-decreasing;}

4) {\it for any} $n \in \mathbb{N}$ {\it the function} $g(n)$ {\it is positive definite and unambiguous;}

5) $\psi \in \Phi$;

6) {\it for some number} $k\in \mathbb{N}_{+}$ {\it the limit}
\begin{gather}\label{eq25}
\lim\limits_{n\rightarrow \infty} \frac{\zeta(n)g(n)}{f(n)\left( \lambda_{k+2}(\psi(n+1))-\lambda_{k+2}(\psi(n))\right) }=S,  \quad S \in \mathbb{R}\setminus\{0\}
\end{gather}
{\it is satisfied; }
\begin{gather}\label{eq26}
7) \quad \lim\limits_{n\rightarrow \infty} \frac{g(n)}{f(n)}=0.
\end{gather}
 {\it Then

1) if for a given number } $k$ {\it satisfying limit \eqref{eq25},} $\exists N:$ $\forall n>N$ {\it the following inequality }
\begin{gather}\label{eq27}
P_{k,n;\psi}:= \left( 1- \left( \frac{a_{n+1}}{a_{n}}\frac{\psi'(n)}{\psi'(n+1) } \prod\limits_{i=0}^{k}\frac{\lambda_{i}(\psi(n+1))}{\lambda_{i}(\psi(n))} \right)^{\frac{1}{\zeta(n)}} \right)\frac{f(n)}{g(n)} \geq R>\frac{1}{S},    \quad R \in \mathbb{R}
\end{gather}
{\it is satisfied, where $R$ is some constant number, then the number series } $\sum\limits_{n=1}^{\infty}a_{n}$ {\it converges;

2) if for a given number } $k$ {\it satisfying limit \eqref{eq25},} $\exists N:$ $\forall n>N$ {\it the following inequality }
\begin{gather}\label{eq28}
\left( 1- \left( \frac{a_{n+1}}{a_{n}} \frac{\psi'(n)}{\psi'(n+1) } \prod\limits_{i=0}^{k}\frac{\lambda_{i}(\psi(n+1))}{\lambda_{i}(\psi(n))} \right)^{\frac{1}{\zeta(n)}} \right)\frac{f(n)}{g(n)} \leq R<\frac{1}{S}
\end{gather}
{\it is satisfied, where $R$ is some constant number, then the number series } $\sum\limits_{n=1}^{\infty}a_{n}$ {\it diverges.}

{\it Theorem 9.1 (Modification test for the case} $ S\rightarrow +\infty ${\it ).}  {\it Let the following conditions be satisfied:}

1) {\it the number sequence } $\{a_{n}\mid n \in \mathbb{N}\}$ {\it is given such that } $a_{n}>0$ {\it for any } $ n \in \mathbb{N}$;

2) {\it for any} $n \in \mathbb{N}$ {\it the function} $f(n)$ {\it is positive definite, unambiguous, increasing and unbounded;}

3) {\it for any} $n \in \mathbb{N}$ {\it the function} $\zeta(n)$ {\it is positive definite, unambiguous and non-decreasing;}

4) {\it for any} $n \in \mathbb{N}$ {\it the function} $g(n)$ {\it is positive definite and unambiguous;}

5) $\psi \in \Phi$;

6) {\it for some number} $k\in \mathbb{N}_{+}$ {\it the limit}
\begin{gather}\label{eq25111}
\lim\limits_{n\rightarrow \infty} \frac{\zeta(n)g(n)}{f(n)\left( \lambda_{k+2}(\psi(n+1))-\lambda_{k+2}(\psi(n))\right) }=+\infty
\end{gather}
{\it is satisfied; }
\begin{gather}\label{eq26111}
7) \quad \lim\limits_{n\rightarrow \infty} \frac{g(n)}{f(n)}=0.
\end{gather}
 {\it Then

1) if for a given number } $k$ {\it satisfying limit \eqref{eq25111},} $\exists N:$ $\forall n>N$ {\it the following inequality }
\begin{gather}\label{eq27111}
P_{k,n;\psi}:= \left( 1- \left( \frac{a_{n+1}}{a_{n}} \frac{\psi'(n)}{\psi'(n+1)} \prod\limits_{i=0}^{k}\frac{\lambda_{i}(\psi(n+1))}{\lambda_{i}(\psi(n))} \right)^{\frac{1}{\zeta(n)}} \right)\frac{f(n)}{g(n)} \geq R>0
\end{gather}
{\it is satisfied, where $R$ is some constant number, then the number series } $\sum\limits_{n=1}^{\infty}a_{n}$ {\it converges;

2) if for a given number } $k$ {\it satisfying limit \eqref{eq25111},} $\exists N:$ $\forall n>N$ {\it the following inequality }
\begin{gather}\label{eq28111}
\left( 1- \left( \frac{a_{n+1}}{a_{n}} \frac{\psi'(n)}{\psi'(n+1) } \prod\limits_{i=0}^{k}\frac{\lambda_{i}(\psi(n+1))}{\lambda_{i}(\psi(n))} \right)^{\frac{1}{\zeta(n)}} \right)\frac{f(n)}{g(n)} <0
\end{gather}
{\it is satisfied, then the number series } $\sum\limits_{n=1}^{\infty}a_{n}$ {\it diverges.}

{\it Theorem 9.2 (Modification test for the case} $ S\rightarrow 0 $ {\it ).} {\it  Let the following conditions be satisfied:}

1) {\it the number sequence } $\{a_{n}\mid n \in \mathbb{N}\}$ {\it is given such that } $a_{n}>0$ {\it for any } $ n \in \mathbb{N}$;

2) {\it for any} $n \in \mathbb{N}$ {\it the function} $f(n)$ {\it is positive definite, unambiguous, increasing and unbounded;}

3) {\it for any} $n \in \mathbb{N}$ {\it the function} $\zeta(n)$ {\it is positive definite, unambiguous and non-decreasing;}

4) {\it for any} $n \in \mathbb{N}$ {\it the function} $g(n)$ {\it is positive definite and unambiguous;}

5) $\psi \in \Phi$;

6) {\it for some number} $k\in \mathbb{N}_{+}$ {\it the limit}
\begin{gather}\label{eq25112}
\lim\limits_{n\rightarrow \infty} \frac{\zeta(n)g(n)}{f(n)\left( \lambda_{k+2}(\psi(n+1))-\lambda_{k+2}(\psi(n))\right) }=0
\end{gather}
{\it is satisfied; }
\begin{gather}\label{eq26112}
7) \quad \lim\limits_{n\rightarrow \infty} \frac{g(n)}{f(n)}=0;
\end{gather}
8) {\it for a given number } $k$ {\it satisfying limit \eqref{eq25112},} $\exists N:$ $\forall n>N$ {\it the following inequality }
\begin{gather}\label{eq28112}
\left( 1- \left( \frac{a_{n+1}}{a_{n}} \frac{\psi'(n)}{\psi'(n+1)} \prod\limits_{i=0}^{k}\frac{\lambda_{i}(\psi(n+1))}{\lambda_{i}(\psi(n))} \right)^{\frac{1}{\zeta(n)}} \right)\frac{f(n)}{g(n)} \leq R
\end{gather}
{\it is satisfied, where $R$ is some constant number.}

{\it Then the number series } $\sum\limits_{n=1}^{\infty}a_{n}$ {\it diverges.}

{\it Theorem 10 (Modification test).} {\it Let the following conditions be satisfied:}

1) {\it the number sequence } $\{a_{n}\mid n \in \mathbb{N}\}$ {\it is given such that } $a_{n}>0$ {\it for any } $ n \in \mathbb{N}$;

2) {\it for any} $n \in \mathbb{N}$ {\it the function} $f(n)$ {\it is positive definite, unambiguous, increasing and unbounded;}

3) {\it for any} $n \in \mathbb{N}$ {\it the function} $g(n)$ {\it is positive definite and unambiguous;}

4) $\psi \in \Phi$;

5) {\it for some number} $k\in \mathbb{N}_{+}$ {\it the limit}
\begin{gather}\label{er9vbn8}
\lim\limits_{n\rightarrow \infty} \frac{g(n)}{f(n)\left( \lambda_{k+2}(\psi(n+1))-\lambda_{k+2}(\psi(n))\right) }=\omega, \quad \omega \in \mathbb{R}\setminus\{0\}
\end{gather}
{\it is satisfied; }
\begin{gather}\label{er10vbn8}
6) \quad \lim\limits_{n\rightarrow \infty} \frac{g(n)}{f(n)}=0.
\end{gather}
 {\it Then

1) if for a given number } $k$ {\it satisfying limit \eqref{er9vbn8},} $\exists N:$ $\forall n>N$ {\it the following inequality }
\begin{gather}\label{er11vbn8}
\tilde{P}_{k,n;\psi}:= \left( 1-  \frac{a_{n+1}}{a_{n}} \frac{\psi'(n)}{\psi'(n+1)} \prod\limits_{i=0}^{k}\frac{\lambda_{i}(\psi(n+1))}{\lambda_{i}(\psi(n))}  \right)\frac{f(n)}{g(n)} \geq \upsilon>\frac{1}{\omega},    \quad \upsilon \in \mathbb{R}
\end{gather}
{\it is satisfied, where $\upsilon$ is some constant number, then the number series } $\sum\limits_{n=1}^{\infty}a_{n}$ {\it converges;

2) if for a given number } $k$ {\it satisfying limit \eqref{er9vbn8},} $\exists N:$ $\forall n>N$ {\it the following inequality }
\begin{gather}\label{er12vbn8}
\left( 1-  \frac{a_{n+1}}{a_{n}} \frac{\psi'(n)}{\psi'(n+1)} \prod\limits_{i=0}^{k}\frac{\lambda_{i}(\psi(n+1))}{\lambda_{i}(\psi(n))}  \right)\frac{f(n)}{g(n)} \leq \upsilon<\frac{1}{\omega}
\end{gather}
{\it is satisfied, where $\upsilon$ is some constant number, then the number series } $\sum\limits_{n=1}^{\infty}a_{n}$ {\it diverges.}

 {\it Theorem 10.1 (Modification test for the case} $ \omega\rightarrow +\infty $ {\it ).} {\it  Let the following conditions be satisfied:}

1) {\it the number sequence } $\{a_{n}\mid n \in \mathbb{N}\}$ {\it is given such that } $a_{n}>0$ {\it for any } $ n \in \mathbb{N}$;

2) {\it for any} $n \in \mathbb{N}$ {\it the function} $f(n)$ {\it is positive definite, unambiguous, increasing and unbounded;}

3) {\it for any} $n \in \mathbb{N}$ {\it the function} $g(n)$ {\it is positive definite and unambiguous;}

4) $\psi \in \Phi$;

5) {\it for some number} $k\in \mathbb{N}_{+}$ {\it the limit}
\begin{gather}\label{er9111vbn8}
\lim\limits_{n\rightarrow \infty} \frac{g(n)}{f(n)\left( \lambda_{k+2}(\psi(n+1))-\lambda_{k+2}(\psi(n))\right) }=+\infty
\end{gather}
{\it is satisfied; }
\begin{gather}\label{er10111vbn8}
6) \quad \lim\limits_{n\rightarrow \infty} \frac{g(n)}{f(n)}=0.
\end{gather}
 {\it Then

1) if for a given number } $k$ {\it satisfying limit \eqref{er9111vbn8},} $\exists N:$ $\forall n>N$ {\it the following inequality }
\begin{gather}\label{er11111vbn8}
\tilde{P}_{k,n;\psi}:= \left( 1-  \frac{a_{n+1}}{a_{n}} \frac{\psi'(n)}{\psi'(n+1) } \prod\limits_{i=0}^{k}\frac{\lambda_{i}(\psi(n+1))}{\lambda_{i}(\psi(n))} \right)\frac{f(n)}{g(n)} \geq \upsilon>0
\end{gather}
{\it is satisfied, where $\upsilon$ is some constant number, then the number series } $\sum\limits_{n=1}^{\infty}a_{n}$ {\it converges;

2) if for a given number } $k$ {\it satisfying limit \eqref{er9111vbn8},} $\exists N:$ $\forall n>N$ {\it the following inequality }
\begin{gather}\label{er12111vbn8}
\left( 1-  \frac{a_{n+1}}{a_{n}} \frac{\psi'(n)}{\psi'(n+1) } \prod\limits_{i=0}^{k}\frac{\lambda_{i}(\psi(n+1))}{\lambda_{i}(\psi(n))}  \right)\frac{f(n)}{g(n)} <0
\end{gather}
{\it is satisfied, then the number series } $\sum\limits_{n=1}^{\infty}a_{n}$ {\it diverges.}

{\it Theorem 10.2 (Modification test for the case} $ \omega\rightarrow 0 $ {\it ).} {\it  Let the following conditions be satisfied:}

1) {\it the number sequence } $\{a_{n}\mid n \in \mathbb{N}\}$ {\it is given such that } $a_{n}>0$ {\it for any } $ n \in \mathbb{N}$;

2) {\it for any} $n \in \mathbb{N}$ {\it the function} $f(n)$ {\it is positive definite, unambiguous, increasing and unbounded;}

3) {\it for any} $n \in \mathbb{N}$ {\it the function} $g(n)$ {\it is positive definite and unambiguous;}

4) $\psi \in \Phi$;

5) {\it for some number} $k\in \mathbb{N}_{+}$ {\it the limit}
\begin{gather}\label{er9112vbn8}
\lim\limits_{n\rightarrow \infty} \frac{g(n)}{f(n)\left( \lambda_{k+2}(\psi(n+1))-\lambda_{k+2}(\psi(n))\right) }=0
\end{gather}
{\it is satisfied; }
\begin{gather}\label{er10112vbn8}
6) \quad \lim\limits_{n\rightarrow \infty} \frac{g(n)}{f(n)}=0;
\end{gather}
7) {\it for a given number } $k$ {\it satisfying limit \eqref{er9112vbn8},} $\exists N:$ $\forall n>N$ {\it the following inequality }
\begin{gather}\label{er12112vbn8}
\left( 1-  \frac{a_{n+1}}{a_{n}} \frac{\psi'(n)}{\psi'(n+1) } \prod\limits_{i=0}^{k}\frac{\lambda_{i}(\psi(n+1))}{\lambda_{i}(\psi(n))}  \right)\frac{f(n)}{g(n)} <\upsilon
\end{gather}
{\it is satisfied, where $\upsilon$ is some constant number.

Then the number series } $\sum\limits_{n=1}^{\infty}a_{n}$ {\it diverges.}

{\it Remark 2.} If the conditions of Theorem 7$^{\ast}$ are met, the following statements are true:

1) Condition \eqref{eq25} in Theorem 9 can be replaced by the condition
$$
\lim\limits_{n\rightarrow \infty} \frac{g(n)\zeta(n)}{f(n) }\frac{\prod\limits_{i=0}^{k+1}\lambda_{i}(\psi(n))}{\psi(n+1)-\psi(n)}=S, \quad S \in \mathbb{R}\setminus\{0\};
$$
2) Condition \eqref{eq25111} in Theorem 9.1 can be replaced by the condition
$$\lim\limits_{n\rightarrow \infty} \frac{g(n)\zeta(n)}{f(n) }\frac{\prod\limits_{i=0}^{k+1}\lambda_{i}(\psi(n))}{\psi(n+1)-\psi(n)}=+\infty;$$
3) Condition \eqref{eq25112} in Theorem 9.2 can be replaced by the condition $\lim\limits_{n\rightarrow \infty} \frac{g(n)\zeta(n)}{f(n) }\frac{\prod\limits_{i=0}^{k+1}\lambda_{i}(\psi(n))}{\psi(n+1)-\psi(n)}=0$;
4) Condition \eqref{er9vbn8} in Theorem 10 can be replaced by the condition
$$
\lim\limits_{n\rightarrow \infty} \frac{g(n)}{f(n) }\frac{\prod\limits_{i=0}^{k+1}\lambda_{i}(\psi(n))}{\psi(n+1)-\psi(n)}=\omega, \quad \omega \in \mathbb{R}\setminus\{0\};
$$
5) Condition \eqref{er9111vbn8} in Theorem 10.1 can be replaced by the condition $\lim\limits_{n\rightarrow \infty} \frac{g(n)}{f(n) }\frac{\prod\limits_{i=0}^{k+1}\lambda_{i}(\psi(n))}{\psi(n+1)-\psi(n)}=+\infty$;
6) Condition \eqref{er9112vbn8} in Theorem 10.2 can be replaced by the condition $\lim\limits_{n\rightarrow \infty} \frac{g(n)}{f(n) }\frac{\prod\limits_{i=0}^{k+1}\lambda_{i}(\psi(n))}{\psi(n+1)-\psi(n)}=0$.

\medskip

\textbf{Conclusion.} The article presents new generalizations of the logarithmic convergence test for number series and the Schlomilch test, various new generalizations of the Jamet test, and modified tests for the convergence of number series. The article consists of two parts. The generalizations of tests for the convergence of number series contained in the first part of the article are simpler compared to the tests in the second part and can most often be used to determine the convergence of series in relatively simple cases.


\begin{thebibliography}{200}
\addcontentsline{toc}{section}{References}
\itemsep=0pt

\bibitem{1}
V. M. Jamet. Sur les series a termes positifs. Nouvelles annales de mathematiques. - 1892. - T.11. - P. 99-103.
\bibitem{2}
Franciszek Prus-Wisniowski. Comparison of Raabe's and Schlomilch's tests, Tatra Mt. Math. Publ. 42 (2009), 119-130.
\bibitem{3}
Prus-Wisniowski, F.: A refinement of Raabe's test, Amer. Math. Monthly 115 (2008), 249-252
\bibitem{4}
Hammond, Christopher NB. "The case for Raabe's test." Mathematics Magazine 93.1 (2020): 36-46.
\bibitem{5}
J.L. Raabe, Note zur Theorie der Convergenz and Divergenz der Reihen. Journal fuur die reine und angewandte Mathematik, 309 - 310, 1834.

\end{thebibliography}
\end{document}